\documentclass[11pt]{article}
\usepackage{amsmath,amsthm,amsfonts,amssymb,amscd}
\usepackage[latin1]{inputenc}

\usepackage{graphicx}

\usepackage{amsmath}

\usepackage{amssymb}

\usepackage{amscd}

\usepackage{bbm}

\newtheorem{Theorem}{Theorem}[section]
\newtheorem{Lemma}[Theorem]{Lemma}

\newtheorem{Proposition}[Theorem]{Proposition}

\newtheorem{prope}[Theorem]{Property}

\newtheorem{Definition}[Theorem]{Definition}

\newtheorem{Remark}[Theorem]{Remark}

\newtheorem{Claim}{Claim}
\newtheorem{Example}{Example}

\def\C{{\mathcal{C}}}
\def\L{{\mathcal{L}}}
\def\G{{\mathcal{G}}}
\def\fa{{\mathcal{F}}}
\def\O{{\mathcal{O}}}

\def\U{{\mathcal{U}}}

\def\po{{\partial}}

\def\Om{{\Omega}}
\def\vr{{\varphi}}
\def\ga{{\gamma}}

\def\la{{\lambda}}

\def\ov{\overline}
\def\al{{\alpha}}

\def\be{{\beta}}

\def\re{{\mathbb{R}}}

\def\bc{{\mathbb{C}}}

\def\tr{\operatorname{{tr}}}
\def\dim{\operatorname{{dim}}}

\def\Hol{\operatorname{{Hol}}}

\def\Diff{\operatorname{{Diff}}}
\def\sing{\operatorname{{sing}}}
\def\Sing{\operatorname{{sing}}}

\title{Codimension one foliations with Bott-Morse singularities I}

\author{Bruno Sc\'ardua and Jos\'e Seade}
\date{}

\begin{document}

\maketitle

\begin{abstract}
We study codimension one (transversally oriented) foliations
$\fa$ on oriented closed manifolds $M$
 having non-empty compact
singular set $\sing(\fa)$  which is locally defined by Bott-Morse
functions. We prove that if the transverse type of $\fa$ at each
singular point  is a center and  $\fa$ has a compact leaf with
finite fundamental group or a component of $\sing(\fa)$
 has codimension $\ge 3$ and finite fundamental group, then all leaves of $\fa$
are compact and diffeomorphic, $\sing(\fa)$ consists of two
connected components, and there is a Bott-Morse function $f:M \to
[0,1]$ such that $f\colon M \setminus \sing(\fa) \to (0,1)$ is a
fiber bundle defining $\fa$ and $\sing(\fa) = f^{-1}(\{0,1\})$.
This yields to a topological description of the type of leaves
that appear in these
 foliations, and also
the type of manifolds admiting such foliations. These results
unify, and generalize, well known results for cohomogeneity one
isometric actions and a theorem of Reeb for foliations with Morse
singularities of center type. In this case each leaf of $\fa$ is a
sphere fiber bundle over each component of $\sing(\fa)$.

\end{abstract}

\renewcommand{\thefootnote}{\fnsymbol{footnote}} 

\section*{Introduction}
\label{section:Introduction}

Cohomogeneity one isometric actions of Lie groups play an
important role in Differential Geometry, particularly in the
Theory of Minimal Submanifolds (see for instance \cite{HL}). A
basic well-known fact about these actions is that whenever the
group and the manifold are compact, if all orbits are principal
then the space of orbits is $S^1$, and if there are special orbits
then there are exactly two of them and the space of orbits is the
interval $[0,1]$. Notice that such an action defines a codimension
one foliation with compact leaves and singular set the special
orbits. Since the action is isometric,  the intersection of the
orbits with a slice $\Sigma$ transverse to a special orbit
corresponds to a Morse singularity of center type.

From the Foliation Theory viewpoint this reminds us of two
important results of Reeb. The first of them is the Complete
Stability Theorem, which states that a transversely oriented
non-singular codimension one foliation having a compact leaf with
finite fundamental group on a closed manifold, is a fibration over
the circle.  The second result concerns foliations with non-empty
singular set. It states that if a codimension one transversely
oriented foliation on a closed manifold has only Morse (isolated)
singularities of center type then there are exactly two such
singularities and the manifold is homeomorphic to a sphere.

In this article we unify these two situations by introducing the
concept of foliation with Bott-Morse singularities. This means
that the singular set of such a foliation is a disjoint union of
compact submanifolds and in a neighborhood of each singular point
the foliation is defined by a Bott-Morse function; so it is a
usual Morse function restricted to each transversal slice. Given
such a foliation, the transverse type of each connected component
of the singular set $\sing(\fa)$ is well-defined, and we can speak
of components of center type, of saddle type, etc., according to
the Morse index of the foliation on a transversal slice.

\vglue.1in

 We prove the following Complete Stability Theorem:

\vskip.2cm

\noindent {\bf Theorem A}. \label{Theorem:CompleteStability}
  {\sl Let $\fa$
be a smooth foliation with Bott-Morse singularities on a closed
oriented manifold $M$ of dimension $m \ge 3$ having only center
type components in $\sing(\fa)$. Assume that  $\fa$ has some
compact leaf $L_o$ with finite fundamental group, or there is a
codimension $\geq 3$ component $N$ with finite fundamental group.
Then  all leaves of $\fa$ are compact, stable, with finite
fundamental group. If, moreover, $\fa$ is transversally
orientable, then $\sing(\fa)$ has exactly two components and
there is a   differentiable Bott-Morse function $f\colon M \to
[0,1]$ whose critical values are $\{0,1\}$ and such that
$f\big|_{M\setminus\sing(\fa)}\colon {M\setminus\sing(\fa)} \to
(0,1)$ is a fiber bundle with fibers the leaves of $\fa$.}

\vskip.2cm

The proof of Theorem A actually  shows that every  compact
transversely oriented foliation with non-empty singular set, all
of Bott-Morse type,  has exactly two components in its
singular set and is given by a Bott-Morse function $f\colon M \to
[0,1]$ as is in the statement.

The first step for proving Theorem A is the following Local
Stability Theorem:

\vglue.1in \vskip.3cm
\noindent {\bf Theorem B}. \label{Theorem:LocalStability} {\sl Let
$\fa$ be a smooth  codimension one foliation on a closed, oriented
manifold $M^m$ having Bott-Morse singularities and let $N^n
\subset \sing(\fa)$ be a {\rm(}compact{\rm)} component with finite
holonomy group {\rm(}e.g., if $N$ has finite fundamental
group{\rm)}. Then there exists a neighborhood $W$ of $N$ in $M$
where $\fa$ is given by a Bott-Morse function $f\colon W \to \re$.
If moreover the transverse type of $\fa$ along $N$ is a center,
then $N$ is stable and the leaves of $\fa$ in $W$, for a suitable
choice of $W$, are fiber bundles over $N$ with fiber $S^{m-n-1}$.
}

\vskip.2cm

Theorem A and its proof lead to the following generalization of
Theorem 1.5 in \cite{LSV}:

\vglue.2in

\noindent{\bf Theorem C}. {\sl Let $\fa$ be a  transversally
oriented, compact foliation  with Bott-Morse singularities on a
closed, oriented, connected manifold $M^m$, $m \ge 3$, with
non-empty  singular set $\sing(\fa)$. Let $L$ be a leaf of  $\fa$.
Then $\sing(\fa)$ has two connected components \,$N_1, N_2$, both
of center type, and one has:
\begin{itemize}

\item[{\rm (i)}]  $M \setminus (N_1 \cup N_2)$ is diffeomorphic to
the cylinder $L \times (0,1)$.

\item[{\rm (ii)}] $L$ is a sphere fiber bundle over both manifolds
$N_1, N_2$ and $M$ is diffeomorphic to the union of the
corresponding disc bundles  over $N_1, N_2$, glued together along
their common boundary $L$ by some diffeomorphism $L \to L$.

\item[{\rm (iii)}] In fact one has a double-fibration
$$N_1 \buildrel{\pi_1}\over{\longleftarrow} L
\buildrel{\pi_1}\over{\longrightarrow} N_2\;,$$
and $M$ is homeomorphic to the corresponding mapping cilynder,
{\it i.e.}, to the quotient space of $(L \times [0,1]) \bigcup
(N_1 \cup N_2) $ by the identifications $(x,0) \sim \pi_1(x)$ and
 $(x,1) \sim \pi_2(x)$.
\end{itemize}}

This yields to a description
 of this type of foliations on manifolds of dimensions $3$ and $4$
 (see Section~\ref{subsection:CompactBottMorse}). In dimension 3 our results 
imply that every foliation as in Theorem C is actually given by a 
cohomogeneity one action of either $SO(3)$ or the 2-torus $S^1\times S^1$.

Unless it is stated otherwise, in this work all manifolds,
bundles, foliations and
 maps are assumed to be of class $C^\infty$. This is just for simplicity, because
essentially everything we say  holds in class $C^r$, for all $r
\ge 1$.

In Section 1 we give the precise definition of foliation with
Bott-Morse singularities and  discuss key-examples of such
foliations.

\vglue.1in This work was done during a visit of the first named
author to the Instituto de Matemáticas of UNAM in Cuernavaca,
Mexico, and a visit of the second named author to IMPA, Rio de
Janeiro, Brazil. The authors want to thank these institutions for
their support and hospitality.

\section{Definitions and examples}
\label{section:Preliminaries}


Let $\fa$ be a codimension one smooth foliation on a manifold $M$
of dimension $m \ge 2$. We denote by $\sing(\fa)$ the singular set
of $\fa$. We say that the singularities of $\fa$ are of {\it
Bott-Morse type\/} if $\sing(\fa)$ is a disjoint union of a finite
number of disjoint compact connected submanifolds, $\sing(\fa) =
\bigcup\limits_{j=1}^t N_j$, each of codimension $\ge 2$,
 and for each $p \in N_j \subset \sing(\fa)$ there exists a neighborhood $V$ of
$p$ in $M$ and a diffeomorphism $\vr\colon V \to P \times D$,
where $P \subset \re^n$ and $D \subset \re^{m-n}$ are discs
centered at the origin, such that $\vr$ takes $\fa|_V$ into the
product foliation $P \times \G$, where $\G=\G(N_j)$ is the
foliation on $D$ given by some Morse function  singularity at the
origin. In other words,  $\sing(\fa) \cap V = N_j \cap V$, \,
$\vr(N_j \cap U) = P\times\{0\} \subset P \times D \subset \re^n
\times \re^{m-n}$ and we can find coordinates $(x,y) =
(x_1,\dots,x_n, y_1,\dots,y_{m-n}) \in V$ such that $N_j \cap V =
\big\{y_1=\cdots=y_{m-n}=0\big\}$ and $\fa|_V$ is given by the
levels of a function $J_{N_j}(x,y) = \sum\limits_{j=1}^{m-n}
\la_j\,y_j^2$ where $\la_j \in \{\pm 1\}$.

The discs $\Sigma_p = \vr^{-1}(x(p)\times D)$ are
transverse to $\fa$ outside $\sing(\fa)$ and the restriction
$\fa|_{\sum_p}$ is an ordinary Morse singularity, whose Morse
index does not depend on the point $p$ in the component $N_j$.
 We shall refer to $\G(N_j) = \fa|_{\sum_p}$ as the {\it
transverse type\/} of $\fa$ along $N_j$. This is a codimension one
foliation in the disc $\Sigma_p$ with an ordinary Morse
singularity at $\{p\}=N_j \cap \Sigma_p$.

If $N_j$ has dimension zero (or if we look at a transversal slice),
then $\fa$ has an ordinary Morse singularity at $p$ and for
suitable local coordinates, $\fa$ is given by the level sets of a
quadratic form $f = f(p)-(y_1^2+\cdots+ y_r^2) + y_{r+1}^2
+\cdots+ y_m^2$\,, where $r \in \{0,\dots,m\}$ is the {\it Morse
index} of $f$ at $p$. The Morse singularity $p$ is a {\it
center\/}  if $r$ is $0$ or $m$, otherwise  $p$ is called a {\it
saddle\/}. In a neighborhood of a center, the leaves of $\fa$ are
diffeomorphic to $(m-1)$-spheres. In a neighborhood of a saddle
$q$, we have conical leaves called {\it separatrices\/} of $\fa$
through $q$, which are
 given by
expressions $y_1^2 +\cdots+ y_r^2 = y_{r+1}^2 +\cdots+ y_m^2\ne0$.
Each such leaf contains $p$ in its closure.

\begin{Definition}
\label{Definition:centersaddlecomponent} {\rm  A component
$N\subset \sing(\fa)$ is of {\it center type} (or just {\it a
center}) if the transverse type $\G(N) = \fa|_{\sum_q}$ of $\fa$
along $N$ is a center. Similarly, the component $N\subset
\sing(\fa)$ is  of  {\it saddle type}  if its transverse type is a
saddle. }

 \end{Definition}

As in the case of isolated singularities, these concepts do not
depend on the choice of orientations. We denote by $C(\fa)\subset
\sing(\fa)$ the union of center type components in $\sing(\fa)$,
and by $S(\fa)$ the corresponding union of saddle components. Of
course  saddles can have different transversal Morse indices;
 this will be relevant for Part II of this article \cite {SS}.

 \begin{Definition}{\rm    We say
that $\fa$ is {\it compact} if every leaf of $\fa$ is compact (and
consequently $S(\fa)=\emptyset$). The foliation $\fa$ is {\it
proper} if every leaf of $\fa$ is closed off $\sing(\fa)$.  }
\end{Definition}

If $\fa$ is proper and  $M$ is compact, then all leaves are
compact except for those containing separatrices of saddles in
$S(\fa)$ and such a leaf is contained in a compact singular
variety $\ov{L}=\ov{L}\cap \sing(\fa) \subset L \cup S(\fa)$.  A
proper foliation on a compact manifold is compact if and only if
$S(\fa)=\emptyset$.

Let $N\subset C(\fa)$ be a component of dimension $k$. Suppose
that the nearby leaves of $\fa$ are compact. We define
$\Omega(N,\fa)=\Omega(N)\subset M$ as the union of $N$ and all the
leaves $L \in \fa$ which are compact and bound a compact invariant
region $R(L,N)$ which is a neighborhood of $N$ in $M$. The region
$R(L,N)$ is equivalent to a fibre bundle with fibre the closed
disc $\ov D^{m-k}$ over $N$, the fibers being transversal to the
leaves of $\fa$.  As we will see, the notion of holonomy of the
singular set, to be introduced in
section~\ref{subsection:Holonomysingularlocalstability},  assures
that if $N$ is of center type and  has finite holonomy group
(e.g., if $\pi_1(N)$ is finite) then $\Omega(N,\fa)$ is an open
subset of $M$.

\begin{Definition} [orientability and transverse orientability]
\label{Definition:orientability} {\rm  Let $\fa$ be a
codimension one foliation with Bott-Morse singularities on $M^m$,
$m \ge 2$. The foliation $\fa$ is {\it orientable\/} if there
exists a one-form $\Om$  on $M$ such that
$\Sing(\fa) = \Sing(\Om)$,  and $\fa$ coincides with the foliation
defined by $\Om=0$ outside the singular set. The choice of such a
one-form $\Om$ is called an {\it orientation\/} for $\fa$. We
shall say that $\fa$ {\it is transversally orientable} if  there
exists a  vector field $X$ on $M$, possibly with
singularities at $\sing(\fa)$, such that  $X$ is transverse to
$\fa$ outside $\sing(\fa)$.}
\end{Definition}

The following basic result is easily proved using
the fact
that we can always choose local orientations for $\fa$, and also orientations
along paths which are null-homotopic.

\begin{Proposition}
\label{Proposition:transversallyorientable} Let $\fa$ be a
codimension one foliation with Bott-Morse singularities on $M^m$,
$m \ge 2$. Suppose $M$ is orientable. Then:

\itemize

\item[{\rm (i)}]  The foliation $\fa$ is orientable if and only if
it is transversally orientable.

\item [{\rm (ii)}] If $M$ is simply-connected, then $\fa$ is
transversally orientable.
\end{Proposition}


\subsection{Examples}

Basic  examples of foliations with Bott-Morse singularities are
given by Bott-Morse functions and by products of Morse foliations
by closed manifolds. Next we give four types of examples of
foliations
 with Bott-Morse singularities.

 \begin{Example}[Fiber bundles]
 \label{Example:fiberbundles} {\rm Let $\widetilde M^{m+k}$ and $ M^m$
be connected oriented manifolds. Let $\fa$ be a foliation with
Bott-Morse singularities on $M$ and let $\pi\colon \widetilde M
\to M$ be a proper submersion. Then the pull-back foliation
$\widetilde \fa=\pi ^* \fa$ has only  Bott-Morse singularities;
hence  $\widetilde \fa$ is a foliation with Bott-Morse
singularities and its transverse type at each component is that of
$\fa$ at the corresponding point.

 For instance, take a vector field on an oriented closed surface $S$ with
non-degenerate singularities, and consider the corresponding foliation $\L$.
Given any $S^1$-bundle  $\pi\colon M \to S$, the pull-back
foliation $\fa=\pi^*(\L)$ has Bott-Morse
singularities on $M$; $\sing(\fa)$ is a union of circles.

 In particular, the  Hopf
fibration $\pi\colon S^3 \to S^2$ gives rise, in this way, to Bott-Morse foliations
on $S^3$. We can consider also $SO(3)$, regarded as the unit
tangent bundle of $S^2$, to get examples on $SO(3) \cong \mathbb R P^3$. }
\end{Example}

\begin{Example}[Mapping cylinders and Lens spaces]
\label{Example:Lensspaces} {\rm \label{Example:Lensspace} Consider
now a closed oriented manifold $L$ that fibers as a  sphere
fiber bundle over two other manifolds $N_1$ and $N_2$, of possibly
different dimensions. Let $E_1, E_2$ be the corresponding disc
bundles. Then each $E_i$ can be foliated by copies of $L$ by
taking concentric spheres in the corresponding fibers. We may now
glue $E_1$ and $E_2$ by some diffeomorphism of the common boundary
$L$ to get a closed oriented manifold $M$ with a foliation with
Bott-Morse singularities at $N_1$ and $N_2$, both of center type.

For instance, take  two solid torii $S^1 \times D^2$, equipped
with the same foliation, and glue their boundaries by a
diffeomorphism that carries a
meridian of the first torus  into a curve on the second which is homologous
 to $q$-meridians and $p$-longitudes, with $p,q \ge 1$ coprime.  We obtain
 foliations with Bott-Morse singularities on the
so-called {\it Lens spaces} $L(p,q)$.
}
\end{Example}

\begin{Example}[Cohomogeneity one actions]\label{Example:projective space}
{\rm As mentioned before, a cohomogeneity
one isometric action leads naturally to compact foliations with
Bott-Morse singularities of center type.

 For instance \cite{LSV}, consider $SO(n+1,\mathbb R)$ as a subgroup of $SO(n+1,\mathbb
C)$. The standard action of this group on $\mathbb C^{n+1}$
defines an action of $SO(n+1,\mathbb R)$ on   $\bc P(n)$, which is
by isometries with respect to the Fubini-Study metric.

 The special orbits are the
complex quadric  $Q_{n-1} \subset \bc P(n)$, of points with
 homogeneous coordinates satisfying $\sum\limits_{j=0}^n z_j ^2=0$,
 and the real projective space $\mathbb R P(n) \subset \bc P(n)$,
consisting of the points which are fixed by the involution in
 $\bc P(n)$ given by complex conjugation.

The principal orbits are copies of the flag manifold
$$F_+^{n+1}(2,1) \cong  SO(n+1,\mathbb
R)/(SO(n-1,\mathbb R)\times (\mathbb Z /2 \mathbb Z))\,,$$
 of oriented $2$-planes in $\mathbb R^{n+1}$
and (unoriented) lines in these planes. Each such orbit splits $\bc P(n)$ in two pieces,
each being a tubular neighborhood of a special orbit.

The case $n =2$ is specially interesting because this provides an equivariant version of the
Arnold-Kuiper-Massey theorem that  $\bc P(2)$ modulo conjugation is the 4-sphere,
see for instance \cite {LSV}. This is also proved in \cite{At-Wh} and \cite {At-Be},
where there are also interesting generalizations of these constructions
and theorem to the quaternionic and the octonian projective planes.
}

\end{Example}

\begin{Example} [Poisson manifolds]
{\rm  A {\it Poisson structure} on a smooth manifold $M$ consists
of a Lie algebra structure on the ring of functions $C^\infty(M)$,
generalizing the classical Poisson bracket on a symplectic
manifold, which satisfies a Leibniz identity in such a way that
$\{\,,h\}$ is a derivation. There is thus a vector bundle morphism
$\psi \colon T^*M\to TM$
 associated with $\{\, , \}$, satisfying an integrability condition,
whose rank at each point is called the rank of the Poisson
structure.

If the rank is  constant, then the integrability
condition implies one has a foliation on $M$, of dimension equal
to the rank, and the tangent space of the foliation is, at each
point $x \in M$, the image of $\psi(T^*_xM)$ in $T_xM$. If the
 rank is not constant, then one still has a generalized foliation
  in the sense of \cite{Su}, {\it i.e.},
 a foliation with singularities at the points where the
rank drops, but at each such point one has a leaf of dimension the
corresponding rank, whose tangent space is again given by
$\psi(T^*_xM)$. The Dolbeault-Weinstein theorem implies that at
such points the transversal structure plays a key role (see \cite {We}).

It would be interesting to study Poisson structures for which the
corresponding foliation has Bott-Morse singularities ({\it cf.} \cite {Du}
 for instance.)}
\end{Example}

\section{Holonomy and local stability}
\label{section:Holonomysingularlocalstability}

 The notion of
stability plays a fundamental role in the classical theory of
(nonsingular) foliations. In what follows we bring this notion
into our framework.

\begin{Definition} {\rm  Let $\fa$ be a (possibly singular) foliation on $M$. A
subset $B \subset M$, invariant by $\fa$, is {\it stable\/} (for
$\fa$) if for any given neighborhood $W$ of $B$ in $M$ there
exists a neighborhood $W' \subset W$ of $B$ in $M$ such that every
leaf of $\fa$ intersecting $W'$ is contained in $W$. }
\end{Definition}

The following technical result comes from the proof of the
Complete Stability theorem of Reeb (cf. \cite{Godbillon}):

\begin{Lemma}
\label{Lemma:reebstability} Let $\fa$ be a codimension one
{\rm(}nonsingular{\rm)} foliation on $M$.
\begin{itemize}
\item[{\rm (i)}] Let $L$ be a compact leaf of $\fa$ and let $L_n$
be a sequence of compact leaves of $\fa$ accumulating to $L$. Then
given a  neighborhood $W$ of $L$ in $M$ one has $L_n\subset W$ for
all $n$ sufficiently  large.

\item[{\rm (ii)}] Denote   by $L_x$ the leaf of $\fa$ containing
$x\in M$ and define $\widehat{M}$ as the subset of points $x\in M$
such that the leaf $L_x$  is compact with  finite fundamental
group.

Then any leaf contained in $\po \widehat {M}$ is closed in $M$.
\end{itemize}
\end{Lemma}

According to \cite{Godbillon}, Proposition 2.20, page 103, in case
$\fa$ is a compact foliation without singularities, stability of a
(compact) leaf is equivalent to finiteness of its holonomy group.
We will extend this result for compact codimension one foliations
with Bott-Morse singularities (see
Proposition~\ref{Proposition:Stabilitycharacterization}) using the
following notion.

\subsection{Holonomy of the singular set}
\label{subsection:Holonomysingularlocalstability}

 Given a component $N \subset
\sing(\fa)$ we consider a collection $\U = \{U_j\}_{j\in J}$ of
open subsets $U_j \subset M$ and charts $\vr_j\colon U_j \to
\vr_j(U_j) \subset \re^m$ with the following properties:

(1) Each $\vr_j\colon U_j \to \vr_j(U_j) \subset \re^m$ defines a
local product trivialization of $\fa$, \, $U_j \cap N$ is a disc
and $\vr_j(U_j)$ is a product of discs.

(2) $\bigcup\limits_{j\in J} U_j$ is an open neighborhood of $N$
in $M$.

(3) If $U_i \cap U_j \ne \emptyset$ then there exists an open
subset $U_{ij} \subset M$ containing $U_i \cup U_j$ and a chart
$\vr_{ij}\colon U_{ij} \to \vr_{ij}(U_{ij}) \subset \re^m$ of $M$,
such that $\vr_{ij}$ defines a product structure for $\fa$ in
$U_{ij}$ and $U_{ij} \cap N \supset (U_i \cup U_j \cap N) \ne
\emptyset$.

\vskip.2cm
\noindent Such a covering $\U$ will be called {\it a chain
adapted\/} to $\fa$ and $N$. When $N$ is compact we can assume
$\U$ to be finite say $\U = \{U_1,\dots,U_{\ell+1}\}$. Suppose now
that $U_j \cap U_{j+1} \ne \emptyset$, \, $\forall\,j \in
\{1,\dots,\ell\}$. In each $U_j$ we choose a transverse disc
$\Sigma_j$\,, \, $\Sigma_j \cap N = \{q_j\}$ such that
$\Sigma_{j+1} \subset U_j \cap U_{j+1}$ if $j \in
\{1,\dots,\ell\}$.  By choice of $\U$ in each $U_j$ the foliation
is given by a smooth function $F_j \colon U_j \to \re$ which is
the natural trivial extension of its restriction to any of the
transverse discs $\Sigma_j$ or $\Sigma_{j+1}$\,.



\par There is a $C^\infty$ local diffeomorphism $\vr_j\colon
(\re,0) \to (\re,0)$ such that $F_{j+1}\big|_{\Sigma_{j+1}} =
\vr_j \circ F_j\big|_{\Sigma_{j+1}}$\,. This implies that
$F_{j+1} = \vr_j\circ F_j$ in $U_j \cap U_{j+1}$ (notice that by
condition (3) if $U_i \cap U_k \ne \emptyset$ then every plaque
$\fa$ in $U_i \backslash N$ intersects at most one plaque of
$U_k\backslash N$).

\begin{Definition}
\label{Definition:Holonomy} {\rm The {\it holonomy map associated
to the chain \/} $\U=\{U_1,\dots,U_\ell\}$ is the local
diffeomorphism $\vr\colon (\re,0) \to (\re,0)$ defined by the
composition $\vr = \vr_\ell \circ\cdots\circ \vr_1$\,. }
\end{Definition}

 Given now a path $c\colon [0,1] \overset{C^0}{\longrightarrow}
N$ we can find a finite chain $\U = \{U_1,\dots,U_{\ell+1}\}$ such
that $\bigcup\limits_{j=1}^{\ell+1} U_j \supset c([0,1])$
 and define {\it the holonomy map of\/} $c\colon [0,1] \to N$
as $\vr = \vr_\ell \circ\cdots\circ \vr_1\colon (\re,0) \to
(\re,0)$. Clearly if $\widetilde c\colon [0,1] \to N$ is
$C^0$-close to $c\colon [0,1] \to N$ and $\widetilde c(0) =
c(0)$,\, $\widetilde c(1) = c(1)$ then $c$ and $\widetilde c$
define the same holonomy map up to isotopy. This shows by a
standard argument that the holonomy map of $c$ is, up to isotopy,
the same holonomy map of any curve $\widetilde c$ homotopic to $c$
in $N$ with $c(0) = \widetilde c(0)$, \, $c(1) = \widetilde c(1)$.
If we now consider closed paths we obtain a map that associates to
each homotopy class $[c] \in \pi_1(N,q_o)$ (where $q_o = c(0)$)
the holonomy map of the path $c\colon [0,1] \to N$. This is indeed
a group homomorphism $\Hol\colon \pi_1(N,q_o) \to
\Diff^\infty(\re,0)$ of the fundamental group of $N$ based at
$q_o$ into the group of germs $C^\infty$ diffeomorphisms fixing
the origin $0 \in \re$. If we move either the base point or the
discs $\Sigma_j$ or else if we change the coverings $\U$ then we
obtain the same homomorphism up to conjugation in
$\Diff^\infty(\re,0)$.

\begin{Definition}
{\rm  We define {\it the holonomy group\/} of the component $N
\subset \sing(\fa)$ as the  quotient of the image of the
homomorphism $\Hol\colon \pi_1(N,q_o) \to \Diff^\infty(\re,0)$ by
conjugation in $\Diff^\infty(\re,0)$.}
\end{Definition}

  In what follows $N \subset
\sing(\fa)$ is compact,  connected and of Bott-Morse type. The
following lemma proves the first statement in Theorem~B.

\begin{Lemma}
\label{Lemma:localstabilitytrivialholonomyfunction} If the
holonomy group of the component $N \subset \sing(\fa)$ is finite
then there is a neighborhood $W$ of $N$ in $M$ where $\fa$ is
given by a smooth function $f\colon W \to \re$.
\end{Lemma}

\begin{proof} We recall (see for instance \cite{Camacho-Lins Neto}
Lemma 5 page 73) that a  finite subgroup of $\Diff^\infty(\re,0)$
is either
 trivial  or has order two and
therefore it is conjugate to the group generated by the involution
$\vr(x)=-x$ in $\Diff(\re,0)$. Assume first that the holonomy is
trivial. The  proof  is by a standard argument of extension by
holonomy. We fix a point $q_o \in N$ and a transverse disc
$\Sigma_{q_o}$ such that $\fa\big|_{\Sigma_{q_o}}$ is given by a
Morse function $f_o\colon \Sigma_{q_o} \to \re$ singular only at
$\{q_o\} = \Sigma_{q_o} \cap N$. Let $q\in N$ be given and
consider a transverse disc $\Sigma_q$ given by a transverse
fibration as in the above definition of holonomy. Fix any curve
$c_q\colon [0,1] \to N$ with $c_q(0)=q_o$ and $c_q(1)=q$. Given a
point $y_o\in \Sigma_{q_o}$ we consider the lift $\tilde
c_{y_o}\colon [0,1]\to L_y$ of the curve $c$ to the leaf $L_{y_o}$
of $\fa$ through the point $y_o$. Put $y=\tilde c_{y_o}(1)\in
\Sigma_q$. We define the value $f(y)=f_o(y_o)$.  By triviality of
the holonomy of $N$ the value $f(y)$ does not depend on the curve
$c_{q}$. Thus we can define a function $f\colon W \to \mathbb R$
in an invariant tubular neighborhood $W$ of $N$ in $M$ with the
following properties:

(i) $f\big|_{\Sigma_{q_o}}=f_o$.

(ii) $f$ is constant along the leaves of $\fa$ in $W$.

(iii) The restriction $f\big|_{\Sigma_q}$ to a transverse disc
$\Sigma_q$ to $N$ at $q\in N$ is conjugate to $f_o$ by a holonomy
map diffeomorphism $h_{c_q}\colon (\Sigma_{q_o},q_o) \to
(\Sigma_q,q)$.

And finally,

(iv)  This extension $f$ is a smooth first integral for $\fa$
which is a submersion in $W\setminus N$.

Assume now that $N$ has holonomy group generated by the real map
$\vr(x)=-x$. Then we can use the same proof of
Lemma~\ref{Lemma:localstabilitytrivialholonomyfunction} above but
replacing $f_o$ by $(f_o)^2$\,. This function $(f_o)^2(x) =
\big(f_o(x)\big)^2$ is invariant by the holonomy $\vr(x) = -x$ and
therefore extends to a well-defined  first integral for $\fa$ in a
neighborhood $W$ of $N$ in $M$.
\end{proof}

\begin{Remark}
\label{Remark:finitegrouponedimension} {\rm In the case the
holonomy has order 2 we cannot assure that the first integral
$f\colon W \to \re$ has connected fibers. Nevertheless, if $\fa$
is transversally oriented then the holonomy of $N$ consists of
orientation preserving elements in $\Diff^\infty(\re,0)$ and
therefore it is finite if and only if it is trivial. This shows
that the order 2 case in the proof of
Lemma~\ref{Lemma:localstabilitytrivialholonomyfunction} does not
occur if $\fa$ is transversally oriented.}
\end{Remark}


\subsection{Proof of  the Local stability}
\label{subsection:Stabilitysingularfoliations}

To prove Theorem~B we use:

\begin{Proposition}
\label{Proposition:Stabilitycharacterization}   Let $\fa$ be a
transversally orientable foliation with Bott-Morse singularities
on $M$. Given a compact component $N \subset \sing(\fa)$ we have:
\begin{itemize}

\item[\rm (i)] If $N$ is of center type and it is a limit of
compact leaves then $N$ is stable.

\item[\rm (ii)] If $\fa$ is compact then $N$ is stable of center
type with trivial holonomy.

\item[\rm (iii)] If $N$ is of center type and the  holonomy group
of  $N$  is finite then $N$ is stable and the nearby leaves are
all compact.

\end{itemize}

\end{Proposition}

\begin{proof}
\noindent First we prove  (i).  Suppose that $N$ is a center and
is a limit of compact leaves of $\fa$ say $N =
\lim\limits_{j\to\infty} L_j$\,. Fix an orientation for $N$.   By
choosing an orientation for $\fa$ in a tubular neighborhood of $N$
on $M$ we may assume that each $L_j$ bounds a region $R_j$ in $M$,
this region is invariant by $\fa$ and therefore, because of the
orientation for $N$, we  can assume that $N \subset R_j$ and $R_{j+1}
\subset R_j$\,, $\forall\, j$ so that $N=\lim\limits_{j\to \infty}
R_j$ as a decreasing limit: the region $R_j$ is invariant.  Since
$\lim L_j = N$ it follows that $\lim R_j = N$ in the Hausdorff
topology.
 Hence every neighborhood $W$ of $N$ in $M$ contains $R_j$ for $j$
big enough; then we can take $W' =$ interior of $R_{j+1}$ so
that $N \subset W' \subset W$ and every leaf $L$ of $\fa$
intersecting $W'$ is contained in $W$. Thus $N$ is stable.

Now we prove (ii).  Assume that $\fa$ is compact. Then, by definition
$N$ is of center type. Since the leaves of $\fa$ are compact,
the holonomy group $\Hol(\fa,N)\subset \Diff^\infty(\mathbb R,0)$
is an orientation preserving  group with finite orbits. This
implies that this group is trivial and therefore $\fa$ has the
product structure in a neighborhood of $N$ in $M$ (cf.
Remark~\ref{Remark:finitegrouponedimension}). Arguing as above we
conclude that given any leaf $L$ close enough to $N$, we may
assume that each $L$ bounds a region $R(L)$ in $M$, this region is
invariant by $\fa$ and such that $\lim_{L\to N} R(L)=N$. Because
of the orientation for $N$  can assume that the above limit is a
decreasing limit so that $N$ is stable.

\vglue.2in \noindent Proof of (iii): As already mentioned, if
$\Hol(\fa,N)$ is finite then it is trivial and $\fa$ has a product
structure in a neighborhood of $N$ which implies that $N$ is
stable with compact nearby leaves.
\end{proof}

\begin{Remark}
\label{Remark:notstable} {\rm For codimension one transversally
oriented {\em nonsingular} foliations, a {\em compact leaf} is
stable if and only if it has trivial holonomy, this is  due to
Reeb~\cite{Reeb5}. Nevertheless, it is not true that a stable
center type component $N\subset \sing(\fa)$ of a foliation with
Bott-Morse singularities $\fa$ necessarily has trivial holonomy. A
counterexample with a one-dimensional component $N\subset
\sing(\fa)$ can be constructed as follows.  Consider the sphere
$S^m$ as the gluing of $S^1\times D^{m-1}$ and $D^2\times S^{m-2}$
through the boundary. On $S^1\times D^{m-1}$ we consider a
non-compact foliation with leaves diffeomorphic to $\mathbb R
\times S^{m-2}$ except for the boundary leaf diffeomorphic to $S^1
\times S^{m-2}$ and on $D^2 \times S^{m-2}$ we consider the
trivial foliation with compact leaves $S^1(r) \times S^{m-2}$.
These foliations glue together through the common boundary leaf
$S^1 \times S^{m-2}$.
 The resulting foliation $\fa$ is  partially depicted in Figure 2
 and has a   non-stable  compact leaf which is diffeomorphic to
$S^1 \times S^{m-1}$.}
\end{Remark}



\begin{proof}
[Proof of  Theorem~B] The first part of the theorem is exactly the
content of
Lemma~\ref{Lemma:localstabilitytrivialholonomyfunction}. Assume
now that  the transverse type of $\fa$ along $N$ is a center. Then
$N$ is stable with compact nearby leaves (cf.
Proposition~\ref{Proposition:Stabilitycharacterization} (iii)). It
remains to prove that these leaves are fibre bundles over $N$ with
fiber $S^{m-n-1}$. The local product structure and the triviality
of the holonomy group give a retraction of a suitable saturated
neighborhood $W$ of $N$ onto $N$  having as fibers transverse
discs $\Sigma$ to $N$. The restriction of this retraction to any
leaf $L\subset W$ gives  a proper smooth submersion of $L$
onto $N$. The fibration theorem of Ehresmann \cite{Godbillon} and
the center type of $N$ give the fibre bundle structure of $L$.
\end{proof}

\section{Complete stability}
\label{section:Completestability}

In this section we  prove the complete stability Theorem A. Our
first step is the following proposition:

\begin{Proposition}
\label{Proposition:completestabilityleaf} Let $\fa$ be a smooth
codimension one  with  Bott-Morse singularities  on a manifold
$M$. Suppose that all components of $\sing(\fa)$ are centers and
there exists a compact leaf $L_o\in \fa$ with finite fundamental
group. Then
 every leaf of $\fa$ is compact with
finite fundamental group. If $\fa$ is transversally orientable
 then all leaves are diffeomorphic to $L_o$.
\end{Proposition}

Using the $2$-fold transversally orientable covering of $\fa$ we
can assume in what follows that $\fa$ is transversally orientable.
To prove Proposition~\ref{Proposition:completestabilityleaf}
denote by $\Omega(\fa)$ the union of leaves $L \in \fa$ which are
compact with finite fundamental group and by $\Omega(L_o)$ the
connected component of $\Omega(\fa)$ that contains the leaf $L_o$.
Let us study the structure of $\Omega(L_o)$. By the Reeb local
stability theorem $\Omega(L_o)$ is open in $M\setminus
\sing(\fa)$. Since $\Omega(L_o)$ is connected and $\fa$ we have
that  all leaves in $\Omega(L_o)$ are diffeomorphic. We claim that
 $\Omega(L_o)=M\setminus \sing(\fa)$,
which obviously implies
Proposition~\ref{Proposition:completestabilityleaf}. This is an
immediate consequence of the following lemma:

\begin{Lemma}
\label{Lemma:completestabilityleaf} We have $\partial \Omega(L_o)
\subset \sing(\fa)$.
\end{Lemma}

\begin{proof} We use  the transverse orientability of $\fa$. If there
is a point $p \in \partial \Omega(L_o) \setminus \sing(\fa)$ then
the leaf $L_p \subset \partial \Omega(L_o)$ is accumulated by
leaves in $\Omega(L_o)$, i.e., by compact leaves with finite
fundamental group.  We claim that  $\ov L_p \cap \sing(\fa) =
\emptyset$. Suppose by contradiction that there exists a component
$N \subset \sing(\fa)$ such that $N \cap \ov L_p \ne \emptyset$.
Then $N\subset \ov L_p \setminus L_p$ and $N$ is a center. We fix
a point $ q\in N$ and a transverse disc $\Sigma \cong \mathbb
R^k$, where $k=\dim M - \dim N$, and such that $\Sigma \cap N =
\{q\}$. We have that $L_p \cap \Sigma$ accumulates at $q$, indeed
$L_p \cap \Sigma$ defines a sequence of $(k-1)$-spheres
$\{S_\nu\}_{\nu \in \mathbb N}$, such that in $\Sigma$ we have
$S_{\nu +1 } \ge S_{\nu}$ with the order given by the inclusion
$B_{\nu +1} \subset B_{\nu}$ and such that $\lim _{\nu \to
+\infty} S_{\nu} =q$. We can assume that $S_{\nu +1} \ne S_\nu$.
Now since $L_p$ is accumulated by leaves in $\Omega(L_o)$ we can
similarly obtain for each fixed $\nu \in \mathbb N$ a sequence of
spheres $\{S_{\nu,j}\}_{j \in \mathbb N}$ that satisfies
$\lim_{j\to \infty} S_{\nu,j} = S_\nu$ and such that $S_{\nu,j}
\subset L_{\nu,j}$ for some leaf $L_{\nu,j}\subset \Omega(L_o)$.
Since the leaf $L_{\nu,j}$ is compact we can assume that
$S_{\nu,j}= L_{\nu,j}\cap \Sigma$ and also we can assume that
$S_\nu < S_{\nu,j_\nu}< S_{\nu-1}$ for some sequence of indexes
$j_1<j_2<...<j_\nu<j_{\nu+1}<...$. Put $L_\nu:= L_{\nu,j_\nu}$.

Now, since $S_{\nu,j_\nu}=L_\nu \cap \Sigma$ and $L_\nu \in
\Omega(L_o)$, it follows that $L_\nu > L_{\nu +1}$ in
$\Omega(L_o)$ and in particular every leaf $L\in \fa$ such that
$L\cap \Sigma$ contains a sphere $S_{\nu,j_\nu}< S_L<S_{\nu+1,
j_{\nu+1}}$ must satisfy $L_{\nu+1}<L<L_\nu$ and in particular
$L\in\Omega(L_o)$. In other words, the leaf $L_p$ cannot satisfy
$L_p\cap \Sigma = \{S_\nu\}$ with $S_{\nu+1,j_{\nu+1}} <
S_\nu<S_{\nu,j_{\nu}}$ and $\lim_{\nu\to \infty} S_\nu =q$,
contradiction. This shows our claim that $\ov L_p \cap \sing(\fa)
= \emptyset$. This claim already shows that $L_p$ must be compact
and since there is a finite covering $L_\nu \to L_p$ for leaves
$L_\nu \in \Omega(L_o)$, it follows that also $L_p$ is compact
with finite fundamental group. Thus $\partial \Omega(L_o) \subset
\sing(\fa)$ completing the proof of
Lemma~\ref{Lemma:completestabilityleaf}.
\end{proof}

The second step in the proof of Theorem~A is:

\begin{Proposition}
\label{Proposition:CompleteStability}  Let $\fa$ be a smooth
codimension one  foliation with Bott-Morse singularities on a
closed manifold $M^m$, $m \ge 3$. Assume that:
\begin{itemize}
\item[\rm (a)] Every component $N \subset \sing(\fa)$ has center
type.  \item[\rm (b)] There is a codimension $\geq 3$ component
$N_o \subset \sing(\fa)$ with  finite fundamental group.
\end{itemize}

\noindent Then $\fa$ is a compact stable foliation. If $\fa$ is
transversally oriented then  its leaves are diffeomorphic.
\end{Proposition}

\begin{proof} Again we can assume that $\fa$ is transversally
orientable. By hypothesis $\sing(\fa) \ne \emptyset$ and contains
a component $N_o^{n_o} \subset \sing(\fa)$ which has a finite
fundamental group and such that (b) is verified. By the Local
Stability Lemma $N_o$ is stable with compact nearby leaves. Also,
since $\fa$ is transversally oriented, the holonomy of $N_o$ is
trivial (see Remark~\ref{Remark:finitegrouponedimension}). This
implies that each leaf of  $\fa$ in $W$ is diffeomorphic to a
$S^{m-n_o-1}$-fibre bundle over $N_o^{n_o}$. By hypothesis (b) we
have $m-n_o \ge 3$. Using the homotopy sequence of the fibration
$S^{m-n_o-1}\hookrightarrow L\to N_o^{n_o}$ we conclude that $L$
has finite fundamental group.  Thus the leaves of $\fa$ in $W$ are
compact with finite fundamental group. Since $\fa$ has codimension
one and is transversally orientable, each such leaf has trivial
holonomy and therefore (by classical Reeb Stability) $\fa$ has a
local product structure in a neighborhood of each such leaf. In
particular, these leaves are diffeomorphic to $L$.

 Let $\O = \{p \in M\backslash\sing(\fa) : L_p$ is
diffeomorphic to $L\}$ then $\O \ne \emptyset$ and $\O$ is an open
subset of $M$ by the Reeb Local Stability Theorem. Indeed, the
proof of the Reeb Complete Stability Theorem shows that every leaf
$L^\prime \subset \po\O$ must accumulate to a  singularity of
$\fa$. But this is not possible because by (a) the components of
$\sing(\fa)$ are stable with compact nearby leaves. Therefore
 $\po\O = \sing(\fa)$ and $M = \O \cup \sing(\fa)$.
This shows that $\fa$ is a compact foliation.
\end{proof}

\begin{Remark}
{\rm  Condition (b) in
Proposition~\ref{Proposition:CompleteStability} is indeed
necessary. For instance,  consider the foliation with Bott-Morse
singularities $\fa$ on $S^2 \times S^2$ given by the product of a
non-periodic flow with exactly two center type singularities on
$S^2$ by the sphere $S^2$, which has non-compact leaves.  }
\end{Remark}

\begin{Proposition}
\label{Proposition:compactfoliationfunction}   Let $\fa$ be a
foliation on a closed connected manifold $M^m$, $m \ge 3$, with
Bott-Morse singularities. Assume that:
\begin{itemize}
\item[\rm (i)] the transverse type of $\fa$ along any component $N
\subset \sing(\fa)$ is a center. \item[\rm (ii)] $\fa$ has some
compact leaf $L_o$ with finite fundamental group.
\end{itemize}

\noindent Then $\fa$ is a compact stable foliation.
\end{Proposition}

\begin{proof} Define the set $\Omega(\fa) \subset M$ as the union of
leaves $L \in \fa$ which are compact and with finite fundamental
group. Then by Reeb Local Stability $\Omega(\fa)$ is open in $M
\backslash \sing(\fa)$. Since $L_o \in \Omega(\fa)$ it follows
that $\Omega(\fa)$ is not-empty and either $\Omega(\fa) = M$ (if
$\fa$ is nonsingular) or $\po\Omega(\fa) \cap \sing(\fa) \ne
\emptyset$. On the other hand, given any component $N \subset
\sing(\fa)$ with $N \cap \po\Omega(\fa) \ne \emptyset$ we must
have $N \subset \po\Omega(\fa)$. Now, since $N \subset
\po\Omega(\fa)$ this implies that $N$ is a limit of compact leaves
of $\fa$ and by Theorem~B the component $N\subset \sing(\fa)$ is
stable. Thus all leaves close enough to $N$ are compact with
finite holonomy group and therefore there is a neighborhood $W$ of
$N$ in $M$ such that $W \backslash N \subset \Omega(\fa)$ and
$\po\Omega(\fa) \cap W = N$. This shows that $\po\Omega(\fa) =
\sing(\fa)$ and therefore $M = \Omega(\fa) \cup \sing(\fa)$.
\end{proof}

 The existence of the function $f\colon M \to [0,1]$ describing $\fa$ in
Theorem~A is a consequence of
Proposition~\ref{Proposition:compactfoliationfunction} to be
proven below. Its proof requires the following lemma.

\begin{Lemma}
\label{Lemma:transversecurve}  Let $\fa$ be a compact foliation on
a closed manifold $M$ having Bott-Morse singularities. Assume that
$\fa$ is transversally orientable  and $\sing(\fa)\ne \emptyset$.
Then $\sing(\fa)$ has exactly two connected components, say $N_1,
\, N_2$, and  there exists an arc $\ga\colon [0,1] \to M$
transverse to $\fa$ such that $\ga(0)\in N_1$, $\ga(1)\in N_2$,
whose image meets every  leaf of $\fa$ at a single point.
\end{Lemma}

\begin{proof} Let us first prove that $\sing(\fa)$ has at most  two
connected components. Take a component $N\subset \sing(\fa)$ and
denote by $\mathcal A(N)$ the subset of $M$ which is the union of
leaves $L\in \fa$ such that $L$ bounds a region $R(L)\subset M$
containing $N$ and such that $\fa\big|_{R(L)\setminus L}$ is a
fibre bundle over $N$. Clearly $N \subset \partial \mathcal A(N)$
and $\partial R(L)\subset \sing(\fa)$. Suppose that there is a
component $N^\prime \subset\partial \mathcal A(N)\setminus N$; let
us prove that $M=N\cup \mathcal A(N)\cup N^\prime$. First we
observe that  $S=N\cup \mathcal A(N)\cup N^\prime$ is an open
subset of $M$. To see this take an invariant  neighborhood $W$ of
$N^\prime$ given by the local stability theorem for $N^\prime$.
Since $N^\prime\subset\partial \mathcal A(N)$ there is a leaf
$L\subset \mathcal A(N)$ which intersects $W$ and therefore is
entirely contained in $W$. By definition of $\mathcal A(N)$ the
leaf $L$ bounds a region $R(L)\subset M$ such that
$\fa\big|_{R(L)\setminus N}$ is a fiber bundle over $N$ and by the
choice of $W$, $L$ bounds a region $R^\prime(L)\subset W$ such
that $\fa\big|_{R(L^\prime)\setminus N^\prime}$ is a fiber bundle
over $N^\prime$. Finally, since $\fa$ has a local product
structure in a neighborhood of $L$ we conclude that $W\setminus
N^\prime\subset \mathcal A(N)$. Thus $S$ is open in $M$. The above
arguments also show that $S$ contains the union of two compact
submanifolds with boundary which are glued along their common
boundary ($L$  above). Hence $S$ equals $M$ and therefore
$\sing(\fa)$ cannot have more components.

To construct the arc $\gamma$ in the statement we first need:

 \begin{Claim}
 \label{Claim:closedcurve} Let $\gamma_o\colon S^1 \to M$ be a closed curve transverse to
 $\fa$ and to $\sing(\fa)$. Then   $\ga_o$ intersects all leaves of $\fa$ and all
components of $\sing(\fa)$.

 \end{Claim}

\begin{proof}[Proof of Claim~\ref{Claim:closedcurve}] Denote by
$\Omega$ the set of all leaves $L\in\fa$ such that $\gamma_o\cap
L\ne \emptyset$. By transversality this is an open set. To see
this set is also closed in $M\setminus \sing(\fa)$ take a
nonsingular point $p\in
\partial \Omega$ and choose an invariant neighborhood $W$ of the leaf
$L_p$ given by the local stability where $\fa$ is trivial. In $W$
any transverse curve to $\fa$ intersects all leaves. Since
$p\in\partial \Omega$ we have $\gamma_o\cap W \ne \emptyset$ and
therefore $\gamma_o\cap L_p \ne \emptyset$. Thus $\Omega$ is
closed in $M\setminus \sing(\fa)$ and $\Omega=M\setminus
\sing(\fa)$. Similar arguments prove that $\gamma_o$ intersects
each component of $\sing(\fa)$.
\end{proof}

Let $X$ be a   vector field transverse to $\fa$ on $M$. Let
$N\subset \sing(\fa)$ be given, we can assume that $X$ is radial
pointing outwards in a neighborhood of $N$. Consider a point $p\in
N\subset \sing(\fa)$ and the orbit $\gamma$ of $X$ whose
$\alpha$-limit is $p$. We consider the $\omega$-limit
$\omega(\gamma)$. Then $\omega(\gamma)$ avoids a neighborhood of
$N$. In fact we have

\begin{Claim}
\label{Claim:omegalimit} $\omega(\gamma)=\{q\}$ where
$q\in\sing(\fa)\setminus N$.
\end{Claim}
\begin{proof}[Proof of Claim~\ref{Claim:omegalimit}]
 Suppose $\omega(\gamma)$ contains some
non-singular point $q$. Then $\gamma$ cuts the leaf $L_q$
infinitely many times. Let us choose two such points $p_1 =
\gamma(t_1)$ and $p_2 = \gamma(t_2)$, \, $t_2 > t_1$\, close
enough to $q$ so that they avoid a neighborhood of $N$ and a path
$\be\colon [0,1] \to L_q$ joining $p_1$ to $p_2$\,. By a classical
argument $\exists\,\delta
> 0$ and a smooth closed curve $\ga_o$\,, transverse to $\fa$ such
that $\gamma$ contains the arc $C = \al([t_1+\delta, t_2-\delta])$
and the complement $\ga_o([0,1])\backslash C$ projects onto $\be$
via a transverse fibration with basis $\be$. Thus we can construct
 a closed curve $\ga_o$ transverse to $\fa$ and which avoids a neighborhood of
 $N$. This contradicts Claim~\ref{Claim:closedcurve}. Therefore
 $\omega(\gamma)\subset \sing(\fa)$. Since $X$ is radial in a neighborhood of $\sing(\fa)$
 we have that $\omega(\gamma)$ has to be a single point, say $\omega(\gamma)=q\in
 \sing(\fa)$. Because $X$ points outwards in a neighborhood of $N$
 we have that $q$ cannot belong to $N$, which proves the claim.
\end{proof}

\begin{figure}[ht]
\begin{center}
\includegraphics[scale=0.6]{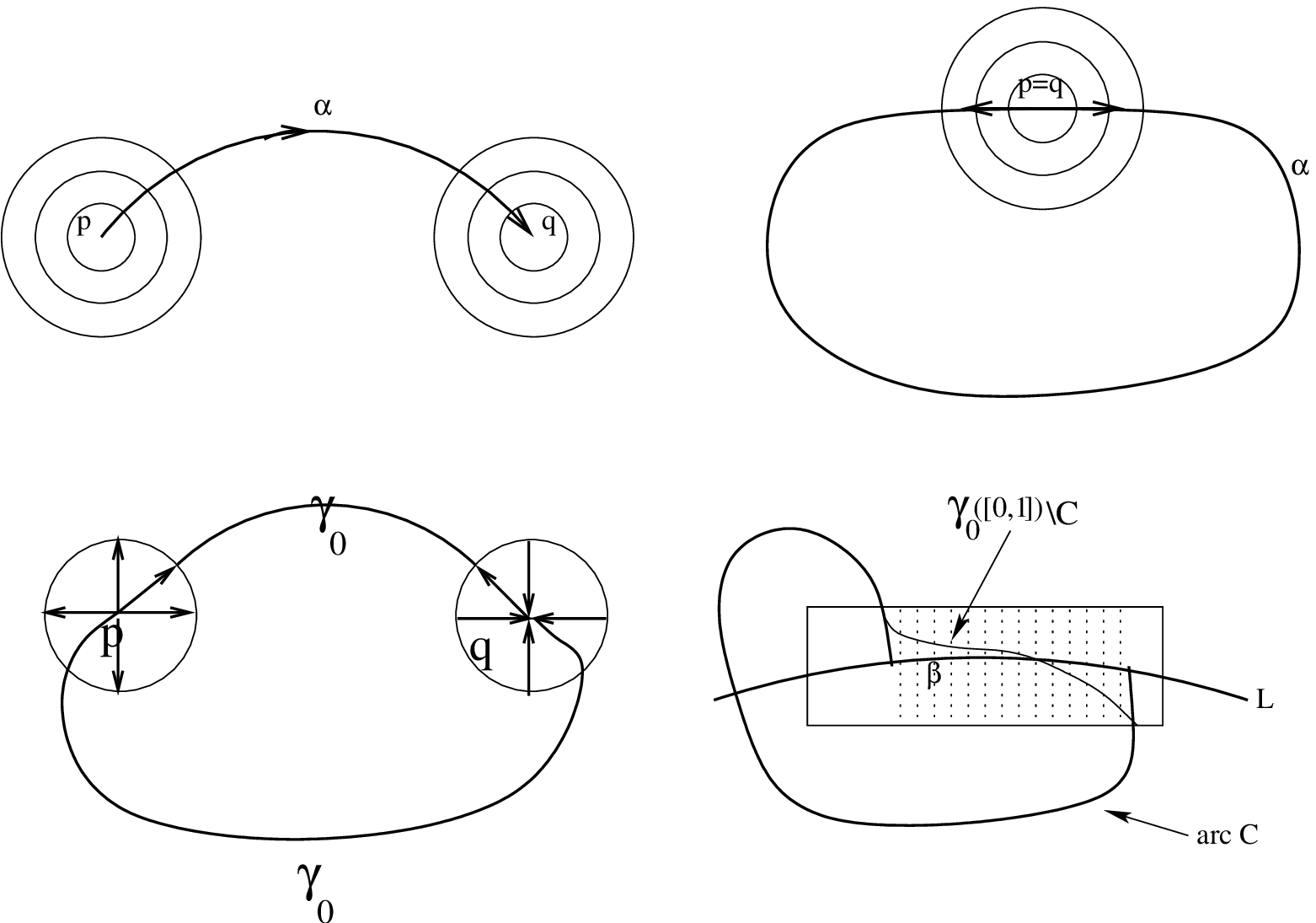}
\end{center}
\caption{}
\end{figure}

Notice that this implies, in particular, that $\sing(\fa)$ has at
least (and therefore exactly) two components $N_1, N_2$. By
Claim~\ref{Claim:omegalimit} we have an arc $\gamma_o\colon
[0,1]\colon M$   such that $\gamma_o(0)\in N_1$ and $\gamma(1)\in
N_2$ and $\gamma(0,1)$ is everywhere transverse to $\fa$. As in
the proof of Claim~\ref{Claim:closedcurve} one can prove that
$\gamma_o$ intersects each leaf of $\fa$. Also by compactness and
transversality $\gamma_o$ intersects each leaf of $\fa$ a finite
number of times.

 If $\ga_o$ cuts a fixed
leaf $L_o$ of $\fa$ a certain number $k \ge 2$ of times then the
proof of Claim~\ref{Claim:omegalimit} above  shows how to modify
$\ga_o$ into a curve $\ga_o'$ that cuts $L_o$ only $k-1$ times.
Therefore we may assume that $\ga_o$ cuts $L_o$ exactly one time.

\begin{Claim}
\label{Claim:onepoint} $\#(\ga_o\cap L) = 1$, \,\, for each leaf
$L$ of $\fa$.
\end{Claim}

\begin{proof}[Proof of Claim~\ref{Claim:onepoint}] Let $\O$ be the set of points $x \in
M\backslash\sing(\fa)$ such that $\#(\ga_o \cap L_x) = 1$.
 By Local Stability we have a local product structure for
$\fa$ around any compact leaf $L$ and therefore $\O$ is open in
$M\backslash\sing(\fa)$. We claim that $\po\O = \sing(\fa)$.
Assume by contradiction that there is a leaf $L \subset \po\O$.
Then again by the local product structure we have $\#(\ga_o \cap
L_x) = \#(\ga_o \cap L)$, $\forall\, x$ close enough to $L$ and
therefore we get a contradiction. This shows that $M = \O \cup
\sing(\fa)$ and proves the claim.
\end{proof}

 Now, by the local
structure of $\fa$ around the singularities we obtain that also
$\#(\ga_o \cap N) = 1$ for any component $N \subset \sing(\fa)$.
This ends the proof of the lemma.
\end{proof}

Theorem~A is now an immediate consequence of Propositions
~\ref{Proposition:CompleteStability} and
~\ref{Proposition:compactfoliationfunction} and
Lemma~\ref{Lemma:transversecurve}.

\subsection{Examples and remarks on stability}
\label{subsection:Examplesapplicationsstability}

  The condition on the codimension of the singular set in
Theorem~A is necessary, for otherwise one can construct examples
of foliations with Bott-Morse singularities with only center type
singularities and non-compact leaves, as for instance the example
in Remark~\ref{Remark:notstable}. A modification of that
construction yields to a foliation with center type singularities which are accumulated
by compact leaves and also by non-compact leaves. For this, let $A^m$ be a
compact  annulus
({\it i.e.}, an $m$-disc minus a smaller $m$-disc in its interior), and
 consider a foliation $\fa_A$  in $S^1\times
A^{m}$ tangent to the boundary
$\po(S^1\times A^{m}) = (S^1 \times S_1^{m-1}) \uplus (S^1\times
S_2^{m-1})$, transverse to the annuli $\{z\} \times A^{m}$, \,
$z \in S^1$ and such that each restriction
$\fa_A\big|_{\{z\}\times A^{m}}$ is equivalent to the trivial
foliation by $(m-1)$-spheres concentric and tangent to the
boundary of $A^{m}$.  We may also choose $\fa_A$ so that each leaf on
$S^1 \times A^{m}$, outside the boundary, is non-compact and
accumulates both components of $\po(S^1\times A^{m})$ as $\fa_o$
above.  Now we consider a sequence of positive numbers $1 = r_1
> r_2 >\cdots> r_j > r_{j+1} >\cdots$ converging to zero. Let
$A_j$ be the annulus of internal radius $r_{j+1}$ and external
radius $r_j$\,. On each solid annulus we put a copy $\fa_{A_j}$ of
$\fa_A$\,. Glue all these foliations in a foliation $\fa_o'$ of
the product $S^1 \times D^{m}$ to get a foliation there, with singular set
$S^1 \times \{pt\}$ of center type. Finally glue two copies
of $\fa_o'$ into a foliation $\fa$ of  $S^1 \times S^m$ with
 two circles $N_1, N_2$ as singular set, both with center types.  Each
component $N_j$ is accumulated by compact leaves (diffeomorphic to
$S^1 \times S^{m-2}$) and also by noncompact leaves (diffeomorphic
to $\re \times S^{m-2}$) as well. In particular, $N_j$ is {\it
not\/} stable and $\fa$ is not compact, although $N_j$ is of center type and
accumulated by compact leaves.

\begin{figure}[ht]
\begin{center}
\includegraphics[scale=0.6]{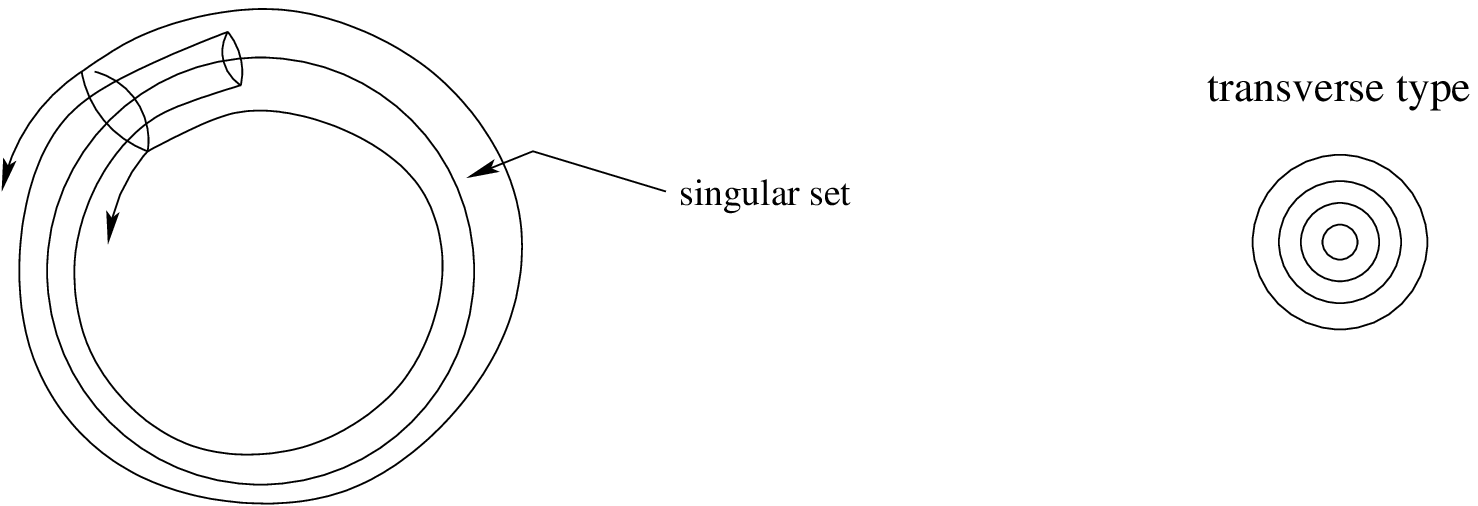}
\caption{}
\end{center}
\end{figure}


\begin{Example}{\rm We decompose $S^3$ as the union of two
solid torii $S^3 = (S_1^1 \times D_1^2) \cup (D_2^2 \times S_2^1)$
with common boundary $S_1^1 \times S_2^1$\,. In the solid torus
$S_1^1 \times D_1^2$ we consider the product foliation $S^1 \times
\C$ where $\C$ is the foliation of the 2-disc by concentric
circles.  We decompose the second solid torus as the union of a
solid torus and a solid annulus with common boundary $S_3^1 \times
S_2^1$, i.e.,  $D_2^2 \times S_2^1 = (A^2\times S_2^1) \cup (D_3^2
\times S_2^1)$. In the solid torus $D_3^2 \times S_2^1$ we put
another trivial foliation $\C \times S^1$. Finally, in the solid
annulus $A^2 \times S_2^1$ we consider a product foliation $\G
\times S_2^1$ where $\G$ is a one-dimensional foliation in $A^2$
as follows:

\begin{figure}[ht]
\begin{center}
\includegraphics[scale=0.4]{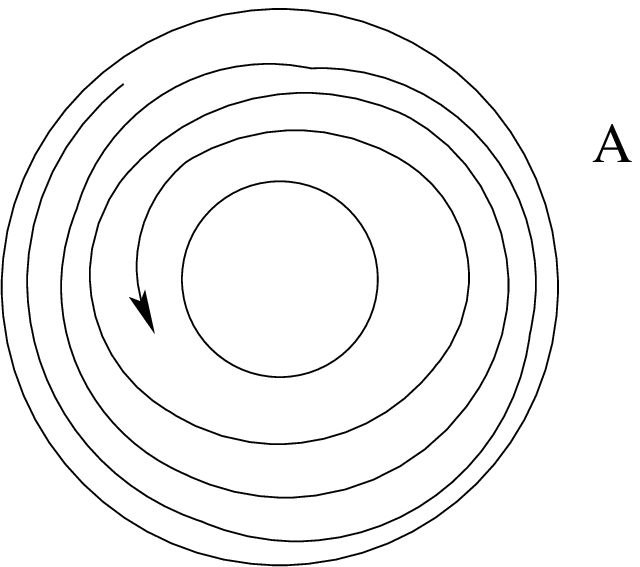}
\end{center}
\caption{}
\end{figure}


\noindent $\G$ has a noncompact leaf accumulating the two circles
in the boundary of $A^2$.  Gluing all together we obtain a
foliation $\fa$ on $S^3$ with singular set $\sing(\fa) =$ union of
two circles which are stable with respect to $\fa$, however $\fa$
is not a foliation by compact leaves due to the noncompact leaves
in $A^2 \times S_2^1$\,. This construction cannot be performed for
dimension $m \ge 4$ as it is implied by
Proposition~\ref{Proposition:compactfoliationfunction}.}\end{Example}

 Let now $\fa$ and $f$ be as in
Theorem~A and assume $m \ge 4$. If $\sing(\fa)$ has some isolated
singularity  then by Reeb complete stability theorem all leaves
are diffeomorphic to $S^{m-1}$. Nevertheless, $\sing(\fa)$ is not
necessarily of dimension zero. For instance, take the classical
Hopf fibration of $S^3$ over $S^2$ with fiber $S^1$. Now consider
the corresponding disc bundle over $S^2$. Its total space $E^4$ is
a four  dimensional manifold with boundary $S^3$.  Using the discs
$D^2$ bounded by the fibers we can construct a foliation with
Bott-Morse singularities $\fa_1$ of $E^4$ having compact leaves
diffeomorphic to $S^3$ and singular set $S^2\cong \mathbb CP(1)$.
Now glue to $E^4$ a four dimensional disc in the obvious way to
obtain the complex projective plane $\mathbb C P(2)$ and a
foliation with Bott-Morse singularities $\fa$ of $\bc P(2)$ with
leaves $S^3$ and singular set $S^2$ union a point. This same
construction generalizes to $\bc P(n)$ regarded as the union of a
$2n$-disc and $\bc P(n-1)$.

\section{Compact foliations with Bott-Morse singularities}
\label{subsection:CompactBottMorse}

Let $\fa$ be a  transversally oriented, compact foliation with
Bott-Morse singularities on the closed, oriented, connected
manifold $M^m$, $m \ge 3$. Notice that
Proposition~\ref{section:Holonomysingularlocalstability}  implies
that  each leaf $L$ of $\fa$ and each component $N$ of the
singular set is stable with  finite holonomy. By
Lemma~\ref{Lemma:transversecurve} one has Theorem C as an immediate
consequence.
 This obviously imposes stringent conditions on both,
the topology of $M$ and $L$. Let us  see what this says when $M$
has dimensions $3$ and $4$. If $m =3$, then $L$ must be a
two-dimensional closed oriented manifold that fibers over another
manifold of dimension 0 or 1, with fiber a sphere. The only
possibilities for $L$ are to be $S^2$, fibered over a point, or
the 2-torus $T = S^1 \times S^1$, since the are no other
$S^1$-bundles over  $S^1$, except for the Klein bottle which is
not orientable. Hence the possibilities for the double-fibration
in Theorem C are:

\begin{itemize}

\item[{\rm (i)}]  If $N_1$ is a point, then $L$ must be a 2-sphere  $S^2$, and this surface
does not fiber over $S^1$, hence $N_2$ must be also a point. This is the classical
case envisaged by Reeb and others,
the leaves are copies of $S^2$ and $M$ is the 3-sphere, regarded as the suspension over $S^2$.

\item[{\rm (ii)}] If $N_1$ is a circle, then $L$ is the torus
$T = S^1 \times S^1$ and $M$
is the result of gluing together two solid torii along their common boundary.
The manifolds one gets in this way are either orientable $S^1$-bundles over $S^2$
(and there is one such bundle
for each integer, being classified by their Euler class), or a lens space
$L(p,q)$, obtained by identifying two solid torii by
  a diffeomorphism of their boundaries that carries a
meridian into a
curve of type $(p,q)$ in $T$.
\end{itemize}

We remark that the hypothesis of having a compact foliation is necessary,
 otherwise the Theorem D does not hold. For instance, decompose $S^3$ as a union
of two solid torii $T_1, T_2$, as usual. Foliate $T_1 = S^1 \times D^2$ by
concentric
torii  $S^1 \times S^1$, and put Reeb's foliation on $T_2$. We get a foliation on
$S^3$
with singular set a circle of center type.

\vskip.2cm

Notice that Theorem C implies:

\vglue.1in \noindent{\bf Theorem D}. {\sl  Let $M$ be a closed
oriented connected $3$-manifold equipped with a transversely
oriented compact   foliation $\fa$ with Bott-Morse singularities.
Then either $\sing(\fa)$ consists of two points, the leaves are
$2$-spheres and $M$ is $S^3$, or $\sing(\fa)$ consists of two
circles, the leaves are torii and $M$ is homeomorphic to a Lens
space or to an $S^1$-bundle over $S^2$. }

\vglue.1in Examples~\ref{Example:fiberbundles} and
\ref{Example:Lensspaces} show that  all $S^1$-bundles over $S^2$
and all Lens spaces admit compact foliations as in Theorem D.

\vglue.1in

\vglue.2in When $m=4$ the list of possibilities for $L$ and $M$ is
larger. For instance, we can foliate $S^4$ in various ways:

\begin{itemize}

\item By $3$-spheres with two isolated centers.

\item By copies of $S^1\times S^2$ with two circles as singular
set.

\item Think of $S^4$ as being the space of real
 $3\times 3$ symmetric matrices $A$ of trace zero and $\tr(A^2) =1$.
  The group $SO(3,\mathbb R)$ acts on
$ S^4$ by $A \mapsto O^t A O$, for a given $O \in SO(3,\mathbb R)$
and $A \in  S^4$. As noticed in \cite{HL} this gives an isometric
action of $SO(3,\mathbb R)$ on the sphere $S^4$ with two
copies of $\mathbb RP(2)$ as singular set. The leaves are copies
of the flag manifold $$F^3(2,1) \cong SO(3,\mathbb R)/(\mathbb Z/2
\mathbb Z \times Z/2\mathbb Z)\cong L(4,1)/(\mathbb Z / 2\mathbb
Z),$$ of (unoriented) planes in $\mathbb R^3$ and lines in these planes.

\item Now consider the complex projective plane $\bc P(2)$.
Thinking of it as being $\bc^2$ union the line
at infinity, one gets a foliation by copies of $S^3$ with an
isolated singularity at the origin and a copy of $S^2 \cong \mathbb C P(1)$ at
infinity.

\item Notice that, as in Example \ref{Example:projective space},
 the group $SO(3,\mathbb R)$ is a subgroup of
$SO(3,\bc)$ and therefore acts on $\bc P(2)$ in the usual way. The
orbits of this action are copies of the Flag manifold $F_+^3(2,1)
\cong SO(3,\mathbb R)/(\mathbb Z/2 \mathbb Z)$, which is a double cover of
$F^3(2,1)$. The
singular set now consists of the quadric $\sum\limits_{j=0}^2 z_j
^2=0$, which is diffeomorphic to $S^2$,  and a copy of $\mathbb R
P(2)$. As in Example \ref{Example:projective space},
this foliation is mapped to the above foliation of $S^4$ by
the projection $\bc P(2)\to \bc P(2)/j \cong
S^4$, where  $j\colon \bc P(2)\to \bc P(2)$ is complex
conjugation (by  \cite{At-Wh}, \cite{At-Be} or  \cite{LSV}).

\end{itemize}

\vglue.1in

 Let us discuss the various possibilities for $L$ and $M$.
Let $N_1$ and $N_2$ be the connected components of $\sing(\fa)$.
If $N_1$ is a point then each leaf $L$ must be $S^3$. We claim
that there are three possibilities for $N_2$: it can be either  a
point,  the $2$-sphere or the projective plane $\mathbb R P(2)$.
Indeed,
 $L$ fibers over $N_2$ with fiber a sphere, and $S^3$ does not fiber over $S^1$.
This implies that $N_2$ has cannot have dimension one. If $N_2$
has dimension two then necessarily  is diffeomorphic to $S^2$ or
to $\mathbb RP(2)$. Thus the possibilities are the following:

\begin{itemize}

\item[{\rm (i.a)}] If $N_2$ is also a point, then $M$ is $S^4$ by Reeb's theorem.

\item[{\rm (i.b)}] If $N_2$ is the 2-sphere then one has a fiber
bundle
$$S^1 \hookrightarrow S^3 \longrightarrow N_2\,;$$
such a bundle necessarily corresponds to a free $S^1$-action on
$S^3$. The effective actions of $S^1$ on 3-manifolds are
classified in \cite{Raymond}, and the only free action on $S^3$ is
the usual one, which yields to the Hopf fibration $S^1
\hookrightarrow S^3 \longrightarrow S^2\,,$ and $M$ is the complex
projective plane $\mathbb C P^2$. Of course the projection $S^2
\to \mathbb RP(2)$ yields to a fibre bundle  $S^1 \hookrightarrow
S^3 \longrightarrow \mathbb RP(2)\,$.

\item[{\rm (ii)}] If $N_1$ is a circle, then $L$ fibers over $S^1
\cong N_1$ with fiber a 2-sphere, so $L$ is $S^1 \times S^2$, and
$N_2$ can be either a circle $S^1$,  $S^2$ or $\mathbb RP(2)$. If
$N_2 \cong S^1$ then both fibrations $L
\buildrel{\pi_i}\over{\longrightarrow} N_i\,,$ $ i =1,2,$
necessarily coincide. Then $M$ is the result of taking two copies
of the corresponding disc bundle, and glued them along their
common boundary $L$ by some diffeomorphism. If  $N_2$ is $S^2$ or
$\mathbb R P(2)$ then  $L$ is a product $S^1 \times S^2$.

\item[{\rm (iii)}] If  $N_1$ and $N_2$ are both  surfaces, then
they can be oriented or not, and $L$ is a closed, oriented Seifert
manifold. The manifolds  $N_1$ and $N_2$ can not be arbitrary,
since $L$ must fiber over both of them simultaneosuly, but there
is  a lot of freedom.
For instance,  notice that we can use the procedure in
Example \ref{Example:Lensspaces} to construct compact
foliations with Bott-Morse singularities whenever we have a
double-fibration as in Theorem C, regardless of whether or not the
hypothesis of Theorem A are satisfied.
\end{itemize}

\bibliographystyle{amsalpha}

\begin{thebibliography}{31}
\frenchspacing

\bibitem{At-Be} M. F. Atiyah, J. Berndt: {\it Projective
varieties, Severi varieties and spheres}; Surveys in Differential
Geometry, Vol. VIII, Boston, MA, 2002, 1--27.



\bibitem{At-Wh} M. F. Atiyah, E.  Witten: {\it  $M$-theory dynamics on a
manifold of $G_2$ holonomy}; Adv. Theor. Math. Phys. 6 (2002), no
1, 1--106.

\bibitem{Bott1} R. Bott: {\it Nondegenerate critical manifolds};
Annals of Math. vol. 60, \# 2, 1954, p. 248-261.


\bibitem{Camacho-Lins Neto} C. Camacho, A. Lins Neto:
Geometric theory of foliations. Translated from the Portuguese by
Sue E. Goodman. Birkhäuser Boston, Inc., Boston, MA, 1985.



\bibitem{Camacho-Scardua} C. Camacho, B.
Sc\'ardua: {\it On codimension one foliations with Morse
singularities}; to appear.


\bibitem{Du} J. P. Dufour: {\it  Lin\'earisation de certaines structures de
Poisson}; J. Diff. Geom. 32 (1990), 415-428.

\bibitem{Godbillon} C. Godbillon: Feuilletages, \'Etudes geom\'etriques,
Birkh\"auser.


\bibitem{Haibao-Rees} D.
Haibao, E. Rees: {\it Functions whose critical set consists of two
connected manifolds}; Bol. Soc. Mat. Mex., vol. 37, 1992, p.
139-149.



 \bibitem{HL}  W. Y. Hsiang,  B. H. Lawson: {\it
 Minimal submanifolds of low cohomogeneity};
 J. Differential Geometry 5 (1971), 1-38.



\bibitem{Ko} A. Kollross: {\it A classification of hyperpolar and
cohomogeneity one actions}; Trans. Amer. Math. 354 (2001),
571--612.



\bibitem{LSV} D. T. L\^e, J. Seade,
A.  Verjovsky: {\it Quadric, orthogonal actions and involuations
in complex projective spaces}; Ens. Math.  49 (2003), no. 1-2,
173-203.


\bibitem{Milnor} J. Milnor: Morse Theory; Ann.
of Math. Studies, Princeton 1968.




\bibitem{Morse1} M. Morse:
{\it The calculus of variations in the large}; American Math. Soc.
Colloquium Publications, 18. AMS, Providence, RI, 1996.


\bibitem{Raymond} F. Raymond: {\it Classification of
the actions of the circle on $3$-manifolds}; Trans. Amer. Math.
Soc. {\bf 131}, 1968, 51--78.


\bibitem{Reeb1}  G. Reeb: {\it
Variétés feuilletées, feuiller voisines}; C.R.A.S. Paris 224,
1947, p. 1613-1614.

\bibitem{Reeb2} G. Reeb: {\it Sur les
points singuliers d'une forme de Pfaff complètement intégrable on
d'une fonction numérique}; C.R.A.S. Paris 222, 1946, p. 847-849.


\bibitem{Reeb4} G. Reeb: {\it Sur certaines propriétés
topologiques des variétés feuilletés}; (Thesis), Publ. Inst. Math.
Univ. Strasbourg 11, p. 5-85, 155-156.

\bibitem{Reeb5} G. Reeb: {\it Sur
certaines propri\'et\'es topologiques des vari\'et\'es
feuillet\'ees};  Publ. Inst. Math. Univ. Strasbourg 11, pp. 5--89,
155--156. Actualit\'es Sci. Ind., no. 1183 Hermann \& Cie., Paris,
1952.

\bibitem{SS} B. Sc\'ardua, J. Seade:
{\it Codimension one foliations with Bott-Morse
singularities II}; in preparation.

\bibitem{Su} H. J.
 Sussmann: {\it Orbits of families of vector fields and integrability
 of distributions};  Trans. Amer. Math. Soc.  180  (1973), 171-188.


\bibitem{We} A. Weinstein: {\it The local structure of Poisson manifolds}; J. Diff. Geom.
18 (1983), 523-557.


\end{thebibliography}

\vglue.2in
\begin{tabular}{ll}
Bruno Sc\'ardua  & \qquad  Jos\'e Seade\\
Instituto de  Matem\'atica & \qquad Instituto de Matem\'aticas \\
Universidade Federal do Rio de Janeiro  & \qquad UNAM\\
Caixa Postal 68530 & \qquad   Unidad Cuernavaca\\
Cidade Universit\'aria & \qquad Av. Universidad s/n,\\
21.945-970 Rio de Janeiro-RJ   & \qquad  Colonia Lomas de Chamilpa\\
BRAZIL &  \qquad C.P. 62210, Cuernavaca, Morelos\\
 \,\, & \qquad M\'EXICO \\
 scardua@impa.br & \qquad jseade@matcuer.unam.mx
\end{tabular}

\end{document}